\documentclass{article}
\usepackage[utf8]{inputenc}
\usepackage{arxiv}
\usepackage{hyperref}       
\usepackage{url}            
\usepackage{amsmath,amssymb,amsthm,graphicx}
\usepackage{booktabs}       
\usepackage{amsfonts}       
\usepackage{nicefrac}       
\usepackage{microtype}      
\usepackage{dsfont}
\usepackage{color}

\newtheoremstyle{thm}
{9pt}
{9pt}
{\itshape}
{}
{\bfseries}
{.}
{ }
{}
\theoremstyle{thm}
\newtheorem{theorem}{Theorem}[section]
\newtheorem{lemma}[theorem]{Lemma}
\newtheorem{corollary}[theorem]{Corollary}

\newtheoremstyle{def}
{9pt}
{9pt}
{}
{}
{\bfseries}
{.}
{ }
{}
\theoremstyle{def}

\newtheorem{remark}[theorem]{Remark}

\newcommand{\R}{\mathbb{R}} 
\newcommand{\E}{\mathbb{E}} 
\newcommand{\PP}{\mathbb{P}} 

    \def\cp{\stackrel{\mathcal{\PP}}{\longrightarrow}}
\renewcommand{\footnoterule}{%
	\kern -3.5pt
	\hrule width \textwidth height 1pt
	\kern 3.5pt
}

\renewcommand{\footnoterule}{%
	\kern -3.5pt
	\hrule width \textwidth height 1pt
	\kern 3.5pt
}

\makeatletter
\def\blfootnote{\xdef\@thefnmark{}\@footnotetext}
\makeatother

\title{Is the Gompertz family a good fit to your data?}
\author{
Dennis Dobler\\
Department of Mathematics\\
Faculty of Science\\
Vrije Universiteit Amsterdam\\
De Boelelaan 1111, NL-1081 HV Amsterdam\\
\href{mailto:d.dobler@vu.nl}{d.dobler@vu.nl}\\
\And
Bruno Ebner\\
Institute of Stochastics \\
Karlsruhe Institute of Technology (KIT) \\
Englerstr. 2, D-76133 Karlsruhe \\
\href{mailto:bruno.ebner@kit.edu}{bruno.ebner@kit.edu}\\
}
\date{\today}

\begin{document}

\maketitle


\begin{abstract}
That data follow a Gompertz distribution is a widely used assumption in diverse fields of applied sciences, e.g., in biology or when analysing survival times. Since misspecified models may lead to false conclusions, assessing the fit of the data to an underlying model is of central importance. We propose a new family of characterisation-based weighted $L^2$-type tests of fit to the family of Gompertz distributions, hence tests for the composite hypothesis when the parameters are unknown. The characterisation is motivated by distributional transforms connected to Stein's method of distributional approximation. We provide the limit null distribution of the test statistics in a Hilbert space setting and, since the limit distribution depends on the unknown parameters, we propose a parametric bootstrap procedure. Consistency of the testing procedure is shown. 
An extensive simulation study as well as applications to real data examples show practical benefits of the procedures: the first data set we analyse consists of lifetimes of fruitflies, the second has been synthetically generated from life tables for women born in Germany in 1948.
\end{abstract}

\section{Introduction}\label{sec:Intro}

The Gompertz distribution was first derived in \cite{G:1825} as a probability model for human mortality. It is widely accepted to describe the distribution of adult lifespans by actuaries and demographers. 
Indeed, as Prof.\ Dr.\ Oliver Ku\ss\ (German Diabetes Center and Heinrich-Heine University D\"usseldorf) pointed out (personal communication): 
``In epidemiology and demography, the Gompertz distribution is used to model the distribution of lifetimes. It is widely accepted that, for ages of 40 years and older, the age at death is nearly perfectly Gompertz-distributed, as seen when comparing empirical (e.g., 5-year) mortality rates with those from a Gompertz fit.''
Beside this prominent application, it is used in fields of applied sciences as biology, see \cite{BM:2016}, and gerontology, see \cite{KMK:1995}, to describe the analysis of survival, in computer science for modeling failure rates, see \cite{HR:2015}, in hydrogen production of energy fuels, see \cite{MZYZ:2006,WYetal:2019} or to describe the walk length of a random self-avoiding walk in the Erd\H{o}s-R\'{e}nyi random graph, see \cite{TBK:2016}.

A first step for serious statistical inference using this model is to assess whether the observed data stems from a distribution being a member of the Gompertz family. Since, even if the assumption of an underlying Gompertz law is true, the true parameters are unknown, it is obvious that this fact has to be incorporated in a test deciding upon the fit of this family of distributions to the data. Hence, in contrast to a simple hypothesis of testing against one fixed Gompertz distribution (including the case that the parameters have been estimated in a first step and then are assumed to be known), a so called composite hypothesis has to be considered. 

The literature for this goodness-of-fit testing problem is hitherto nearly non existent. The only exception is the article \cite{AM:2016}. The authors provide a comparative simulation study for the classical Anderson-Darling test and some extension due to \cite{SSA:1990}, as well as a correlation coefficient type test and a nested test against the truncated generalised extreme value distribution for the minimum. The cited article does not consider testing procedures for the composite case, since the authors first fit the parameters to a Gompertz law and then perform the simple hypothesis tests. No asymptotic theory for the testing procedures is derived. So we conclude that the composite case has hitherto not been treated in the literature. 

In this article, we present new tests of fit to the Gompertz family of distributions based on a characterisation of the Gompertz law. The idea of using characterisations for deriving goodness-of-fit goes back to \cite{L:1953} and is the basis for powerful procedures; for details see \cite{N:2017}. The aim of this article is to propose the first characterisation based test of fit to the composite hypothesis that the data stems from (any) Gompertz law and provide asymptotic theory under the null hypothesis. Since the limit null distribution depends on the unknown shape parameter a parametric bootstrap procedure is presented and we provide the first comparative Monte Carlo simulation study for this setting.
To be precise, let $GO(\eta,b)$ denote the Gompertz distribution defined by the probability density function
\begin{equation}\label{eq:density}
f(x;\eta,b)=b\eta\exp(\eta+bx-\eta e^{bx}),\quad x\ge0,    
\end{equation}
where the rate $b>0$ is a scale parameter and $\eta>0$ is a shape parameter. The cumulative distribution function is given by
\begin{equation*}
F(x;\eta,b)=1-\exp(-\eta(e^{bx}-1)),\quad x\ge0,  
\end{equation*}
and $F(x;\eta,b)=0$ otherwise. We write $\mathbf{\mbox{GO}}=\{GO(\eta,b):\eta,b>0\}$ for the family of Gompertz distributions, see \cite{MO:2007}, Chapter 10, for details on the parametric family $\mathbf{\mbox{GO}}$. Let $X,X_1,X_2,\ldots$ be positive,  independent and identically distributed (i.i.d.) random variables defined on a common probability space $(\Omega,\mathcal{A},\mathbb{P})$, and denote the distribution of $X$ by $\mathbb{P}^X$. We test the composite hypothesis
\begin{equation}\label{eq:H0}
H_0:\;\mathbb{P}^X\in\mathbf{\mbox{GO}}
\end{equation}
against general alternatives.


This paper is organised as follows. In Section~\ref{sec:prep}, we introduce the family of Gompertz distributions, provide a Stein characterisation of these distributions, and propose a goodness-of-fit test statistic $T_n$ that is based on the Stein characterisation. 
In Section~\ref{sec:lnd}, we analyse the distribution of $T_n$ under the null hypothesis and develop a parametric bootstrap scheme for approximating this null distribution.
In particular, we derive the asymptotic distribution of $T_n$ under the assumption of a Gompertz law using a Hilbert space framework, and we show that the parametric bootstrap procedure is well calibrated.
This is followed by Section~\ref{sec:cons} where we consider the behaviour of the test statistic under alternatives: the procedure is consistent under a convergence assumption for the estimators.
We present a competitive Monte Carlo simulation study under both the null and alternative hypothesis in Section~\ref{sec:simu} and compare the new tests to the classical empirical distribution based methods.
The results indicate that the new test is a strong competitor to the classical procedures. The procedure is also applied to real data sets from biology and demography in Section~\ref{sec:real}.
We conclude the paper by reflecting our findings and stating open problems in Section \ref{sec:CaO}.
All proofs and lengthy derivations are offered in the appendices.

\section{Characterisation of the Gompertz law and the new test statistic}
\label{sec:prep}

This article studies a test procedure for (\ref{eq:H0}) based on a characterisation of the family of Gompertz distributions due to \cite{BE:2021}, Corollary 3. This type of characterisation is related to distributional characterisations in Stein's method (for an introduction to the topic we refer to \cite{CGS:2011}) and the so-called density approach, see \cite{LS:2013,SDHR:2004}. For the sake of completeness, we state the characterisation, a proof is found in Appendix \ref{app:PT1}.
\begin{theorem}\label{thm:char}
Let $X$ be a positive random variable with cumulative distribution function $G$ and $\mathbb{E}|X| < \infty$. Define for $\eta, b > 0$ the function $T^X: \mathbb{R}\rightarrow \mathbb{R}$ 
\begin{equation*}
T^X (s) = \left\{\begin{array}{cc}
     \mathbb{E}\big[(\eta b e^{bX} - b)\min\{X, s\}\big], & s > 0,  \\
     0, & s \leq 0.
\end{array}\right. 
\end{equation*}
Then $X\sim GO(\eta, b)$, if and only if $T^X \equiv G$ on $\mathbb{R}$.
\end{theorem}

In the following, we write $\cp$ for convergence in probability, and $(\widehat{\eta}_n,\widehat{b}_n)$ for consistent estimators of $(\eta,b)$. We assume throughout that, for all considered distributions, we have $(\widehat{\eta}_n,\widehat{b}_n)\cp(\eta_0,b_0)\in(0,\infty)^2$ as $n\rightarrow\infty$, i.e. if $X\sim GO(\eta,b)$, we have $(\widehat{\eta}_n,\widehat{b}_n)\cp(\eta,b)$ as $n\rightarrow\infty$. In view of the scale invariance of the Gompertz family, we set $Y_{n,j}=\widehat{b}_nX_j$, $j=1,\ldots,n$ and assume that $\widehat{b}_n$ is a scale equivariant estimator of $b$, i.e., that we have
\begin{equation*}
    \widehat{b}_n(\beta X_1,\ldots,\beta X_n)=\widehat{b}_n(X_1,\ldots,X_n)/\beta
\end{equation*}
and that $\widehat{\eta}_n$ is scale invariant, i.e. that 
\begin{equation*}
    \widehat{\eta}_n(\beta X_1,\ldots,\beta X_n)=\widehat{\eta}_n(X_1,\ldots,X_n)
\end{equation*}
holds for all $\beta>0$. This implies $\widehat{b}_n(Y_{n,1},\ldots,Y_{n,n})=1$, such that under $H_0$ the distribution of the random variables $Y_{n,j}$ should be close to a $GO(\eta,1)$ distribution. Hence, based on Theorem \ref{thm:char} we propose the weighted $L^2$-type statistic
\begin{equation*}
    T_n = n\int_0^\infty V_n^2(s)\, w(s) \;\mbox{d}s,
\end{equation*}
where
\begin{align} \label{eq_V}
    V_n (s) = \frac{1}{n} \sum_{j = 1}^n (\widehat{\eta}_n e^{Y_{n, j}} - 1)\min\{Y_{n, j}, s\} - \frac{1}{n} \sum_{j = 1}^n \mathbf{1}\{Y_{n, j} \leq s\}, \quad s > 0,
\end{align}
and $w(\cdot)$ is a continuous positive weight function, with 
\begin{align}\label{eq:w1}
    \int_0^\infty (s^2 +1)w (s) \,ds < \infty,
\end{align}
and
\begin{equation}\label{eq:w2}
    n \int_0^\infty \big|w(\widehat{b}_n s) - w(s)\big|^3 w^{-2} (s) \,ds \stackrel{\mathbb{P}}{\longrightarrow} 0, \qquad n \rightarrow \infty.
\end{equation}
We reject the hypothesis $H_0$ for large values of $T_n$. Note that $T_n$ only depends on the rescaled data $Y_{n,j}$, $j=1,\ldots,n$, and as a consequence it is invariant due to scale transformations of the data, i.e. w.r.t. transformations of the form $x\mapsto \beta x$, $\beta>0$.

It is straightforward to show that for all $a>0$ the weight function $w_a(s)=\exp(-as)$, $s>0$, satisfies \eqref{eq:w1} and \eqref{eq:w2}. Direct calculations lead to the numerical stable integration free representation 
\begin{align*}
    T_{n, a} = &n\int_0^\infty V_n^2 (s) w_a(s) \,ds\\
    = &\frac2n \sum_{0 \leq i < j \leq n} \Bigg[ \frac{G_j}{a^3} e^{-aY_{(i)}} \big(-aG_iY_{(i)} - 2G_i -a^2Y_{(i)} - a\big) \\
    &+ \frac{1}{a^2} e^{-aY_{(j)}} \big( -G_iG_jY_{(i)} - aG_iY_{(i)} + G_j + a \big) + 2G_iG_j\Bigg] \\
    &+ \frac1n \sum_{j=1}^n \Bigg[ \frac{1}{a^3} e^{-aY_{(j)}} \big(-2aG_j^2Y_{(j)} - 2G_j^2 - 2a^2G_jY_{(j)} + a^2\big) + 2G_j^2 \Bigg],
\end{align*}
where $Y_{(j)}$ stands for the $j$th order statistic of $Y_{n, 1}, \ldots, Y_{n, n}$ and $G_j = \widehat{\eta}_n e^{Y_{(j)}}-1$, $j = 1, \ldots, n$. Note that the proposed test is in the spirit of Stein goodness-of-fit tests, see  \cite{SReview:2023} and the references therein for details on the general approach which is also applicable for other families of distributions. In \cite{BEN:2022} discrete analogs to Theorem \ref{thm:char} are derived and applied to testing the fit to families of discrete distributions. 

\section{Limit null distribution and bootstrap procedure}\label{sec:lnd}
In this section, we derive the asymptotic distribution under the null hypothesis. Due to the $L^2$-structure of the test statistic, a convenient setting is the separable Hilbert space $\mathbb{H}=L^2([0,\infty),\mathcal{B},w(t){\rm d}t)$ of (equivalence classes of) measurable functions $f:[0,\infty) \rightarrow \R$
satisfying $\int_0^\infty |f(t)|^2 \,w(t)\, {\rm d}t < \infty$. Here, $\mathcal{B}$ denotes the Borel sigma-field on $[0,\infty)$. The scalar product and the norm in $\mathbb{H}$ will be denoted by
\begin{equation*}
\langle f,g \rangle_{\mathbb{H}} = \int_0^\infty f(t)g(t)\,w(t)\, {\rm d}t, \quad \|f\|_{\mathbb{H}} = \langle f,f \rangle_{\mathbb{H}}^{1/2}, \quad f,g \in \mathbb{H},
\end{equation*}
respectively. In view of the scale invariance of $T_n$ we assume in the following that  $X_{n,1},\ldots,X_{n,n}$ is a triangular array of rowwise i.i.d. random variables, and suppose $X_{n,1}\sim GO(\eta_n,1)$ for a sequence of positive parameters $(\eta_n)$, where $\lim_{n\rightarrow\infty}\eta_n=\eta>0$. In the following, we assume that the estimators $(\widehat{\eta}_n,\widehat{b}_n)$ allow linear representations 
\begin{align}
    \sqrt{n} (\widehat{\eta}_n - \eta_n) &= \frac{1}{\sqrt{n}} \sum_{j = 1}^n \psi_1 (X_{n, j}, \eta_n) + o_{\mathbb{P}}(1),\label{eq:psi_11}\\
    \sqrt{n} (\widehat{b}_n - 1) &= \frac{1}{\sqrt{n}} \sum_{j = 1}^n \psi_2 (X_{n, j}, \eta_n) + o_{\mathbb{P}}(1),\label{eq:psi_21}
\end{align}
where $\psi_1$ und $\psi_2$ are measurable functions with
\begin{align}
    &\mathbb{E}[\psi_1 (X_{n, 1}, \eta_n)] = 0,  &&\mathbb{E}[\psi_2 (X_{n, 1}, \eta_n)] = 0, \label{eq:psi1}\\
    &\mathbb{E}[\psi_1^2 (X_{n, 1}, \eta_n)] < \infty,  && \mathbb{E}[\psi_2^2 (X_{n, 1}, \eta_n)] < \infty, \label{eq:psi2}
\end{align}
and
\begin{align} \label{eq:psi3}
    \lim_{n \rightarrow \infty} \mathbb{E}[\psi_1^2(X_{n, 1}, \eta_n)] = \mathbb{E}[\psi_1^2(X, \eta)], \qquad \lim_{n \rightarrow \infty} \mathbb{E}[\psi_2^2(X_{n, 1}, \eta_n)] = \mathbb{E}[\psi_2^2(X, \eta)].
\end{align}
Here and in the following $o_{\mathbb{P}}(1)$ stands for a term that converges to 0 in probability. An example for $(\widehat{\eta}_n,\widehat{b}_n)$ satisfying these assumptions are the maximum-likelihood estimators. This fact can be proved by straightforward calculations showing the existence of the inverse Fisher information matrix and applying Theorem 5.39 in \cite{V:1998}.

\begin{theorem}\label{thm:H0vert}
    Under the triangular array $X_{n, 1}, \ldots, X_{n, n}$, we have  
    \begin{equation*}
    T_n = n\Vert V_n\Vert_{\mathbb{H}}^2 \stackrel{\mathcal{D}}{\longrightarrow} \Vert \mathcal{W}\Vert_{\mathbb{H}}^2, \quad n \rightarrow \infty,
    \end{equation*}
    where $\mathcal{W}$ is a centred Gaussian random element in $\mathbb{H}$ having covariance kernel
    \begin{align*}
    \mathcal{K}_\eta (s, t) = &(s-2) F(s; \eta,1) + s(1-t)F(t; \eta,1) + 2\mathbb{E}\big[X\mathbf{1}\{X \leq s\}\big]+ \mathbb{E}\big[X^2\mathbf{1}\{X \leq s\}\big] \\
    & + \mathbb{E}\big[X\mathbf{1}\{X \leq t\}\big]+ \mathbb{E}\big[R (s; \eta) r(t; \eta)\big] + \mathbb{E}\big[R (t; \eta) r(s; \eta)\big] + \mathbb{E}\big[r (s; \eta) r(t; \eta)\big]
    \end{align*}
    for $0 < s \leq t < \infty$, where $X\sim GO(\eta,1)$ and 
    \begin{align*}
        R (s; \eta) = &(\eta X e^{X} - X - 1) \mathbf{1}\{X \leq s\} + s(\eta e^{X} - 1)\mathbf{1}\{X > s\}, \\
    r (s; \eta) = &\psi_1(X, \eta) \frac{1}{\eta}\big(1 + \mathbb{E} \big[X \mathbf{1}\{X \leq s\}\big]\big)+ \psi_2(X, \eta) \big(1 + \mathbb{E}\big[X\big] + \mathbb{E}\big[X\mathbf{1}\{X \leq s\}\big] + \mathbb{E}\big[X^2\mathbf{1}\{X \leq s\}\big]\big).
    \end{align*}
\end{theorem}

The distribution of $\|\mathcal{W}\|_{\mathbb{H}}^2$ is known to have the equivalent representation $\sum_{j=1}^\infty\lambda_j(\eta)N_j^2$, where $N_1,N_2,\ldots$ are independent, standard normally distributed random variables, and $\lambda_1(\eta),\lambda_2(\eta),\ldots$ is the decreasing series of non-zero eigenvalues of the integral operator
\begin{equation*}
K:\mathbb{H}\rightarrow\mathbb{H},\quad f\mapsto Kf(\cdot)=\int_0^\infty \mathcal{K}_\eta(\cdot,t)f(t)\,w(t)\,\mbox{d}t.
\end{equation*}
This operator obviously depends on the true unknown parameter $\eta>0$. To calculate the eigenvalues $\lambda$ of $\mathcal{K}$, one has to solve the homogeneous Fredholm integral equation of the second kind
\begin{equation}\label{eq:inteq}
\int_0^\infty \mathcal{K}_\eta(x,t)f(t)\,w(t)\,\mbox{d}t=\lambda f(x),\quad x>0,
\end{equation}
see, e.g., \cite{KS:1947}. Due to the complexity of the covariance kernel, it seems hopeless to find explicit solutions of (\ref{eq:inteq}) and hence formulae for the eigenvalues. Although numerical and stochastic approaches to approximate the eigenvalues can be found in the literature, see for example \cite{BMNO:2020,EH:2021,S:1976}, the complexity of $\mathcal{K}_\eta$ indicates that even these approaches are hard to apply. Furthermore, since the true parameter $\eta$ is unknown in practice, the limiting null distribution cannot be used to derive critical values of the test. A solution to this problem is provided by a parametric bootstrap procedure as suggested in \cite{H:1996} and which is stated as follows:
\begin{enumerate}
  \item[(1)] Compute $\widehat{\eta}_n=\widehat{\eta}_n(X_1,\ldots,X_n)$. 
  \item[(2)] Conditionally on $\widehat{\eta}_n$ and $\widehat{b}_n=1$, simulate $B$ bootstrap samples $X_{j,1}^*,\ldots,X_{j,n}^*$, i.i.d. from $GO\left(\widehat{\eta}_n,1\right)$, and compute $T_{n,j}^*=T_n(X_{j,1}^*,\ldots,X_{j,n}^*)$, $j=1,\ldots,B$.
  \item[(3)] Denote by 
  \begin{equation*}
  H_{n, B}^* (s) = \frac{1}{B} \sum_{j = 1}^B \mathbf{1} \{T_{n, j}^* \leq s\}, \qquad s > 0,
  \end{equation*}
  the empirical distribution function of $T_{n,1}^*,\ldots,T_{n,B}^*$ and derive the empirical $(1-\alpha)$-quantile $c_{n,B}^*(\alpha)$.
  \item[(4)] Reject the hypothesis (\ref{eq:H0}) at level $\alpha$ if $T_n(X_1,\ldots,X_n)>c_{n,B}^*(\alpha)$.
\end{enumerate}
Note that for each computation of $T_{n,j}^*$, parameter estimation has to be done separately for each $j$ and clearly the tests depend on the rescaled bootstrap data.
The following theorem gives the final justification for the right-tailed test procedure for testing $H_0$ based on the test statistic $T_n$.
\begin{theorem} \label{satz:boot1}
Let $X \sim GO(\eta_n, 1)$, where  $(\eta_n)$ is again a positive sequence with $\lim_{n\rightarrow\infty} \eta_n=\eta > 0$ and $X_1, \ldots, X_n$ be i.i.d copies of $X$. Further denote by $H_{n, B}^*$ and $c_{n, B}^*$ the quantities from the bootstrap procedure. Then, we have
\begin{equation*}
\mathbb{P}_\eta \big(T_n > c_{n, B}^* (\alpha)\big) \longrightarrow \alpha, \qquad\text{as } n, B \rightarrow \infty.
\end{equation*}
\end{theorem}
\section{Consistency}\label{sec:cons}
In this section, we assume that $X$ is a positive, non-degenerate random variable with an arbitrary distribution such that $\E[X]<\infty$ and $\E[X\exp(X)]<\infty$. In view of the scale invariance of $T_n$, we assume $b_0=1$, i.e. that $(\widehat{\eta}_n,\widehat{b}_n)\cp(\eta_0,1)$, $\eta_0>0$, for $n\rightarrow\infty$.
\begin{remark}
Note that, under certain regularity conditions, maximum likelihood estimators in misspecified models are known to be consistent for the minimiser of a Kullback-Leibler information criterion; see Theorem~2.2 in \cite{W:1982}.
In the following, we will assume that such a consistency holds.
\end{remark}
\begin{theorem} \label{thm:cons}
Let $X_1, \ldots, X_n$ be i.i.d. copies of $X$. Then
\begin{align*}
    \frac{T_n}{n} \stackrel{\mathbb{P}}{\longrightarrow} \Delta_{\eta_0},\quad n\rightarrow\infty,
\end{align*}
where $\Delta_{\eta_0} = \| \Delta_{\eta_0}^* \|_{\mathbb{H}}^2$ and $\Delta_{\eta_0}^* (s) = \mathbb{E}\big[ (\eta_0 e^X - 1) \min \{X, s\}\big] - \mathbb{P} (X \leq s),\quad s>0$.
\end{theorem}
The following corollary shows that the parametric bootstrap testing procedure based on $T_n$ is consistent against a broad class of alternatives. 
\begin{corollary}\label{cor:cons}
Under the standing assumptions with the notations from Section \ref{sec:lnd}, we have if $\mathbb{P}^X\not\in\mathbf{\mbox{GO}}$ 
\begin{equation*}
    \mathbb{P} \big(T_n > c_{n, B}^* (\alpha)\big) = 1, \qquad n, B \rightarrow \infty.
\end{equation*}
\end{corollary}

\section{Simulation study}\label{sec:simu}
We assess the practical usefulness of the new test with the help of an extensive simulation study.
We chose the significance level $\alpha=5\%$, sample sizes $n=20,50,100$, and repeated each test for each simulation scenario 10,000 times, where each test was based on 2,000 parametric bootstrap iterations. 
The test was applied for several choices of the tuning parameter, $a \in \{0.1,0.25,0.5,0.75,1,1.5,2,3,5,10\}$.
Additionally, we considered the following competitor goodness-of-fit tests, which are also scale-invariant because they are based on the empirical distribution function $\widehat F_n$ of the re-scaled data, $Y_j = \widehat b_n X_j$:
\begin{align*}
    & \text{(Kolmogorov-Smirnov)}  &&  \mbox{KS} = \sup_x |\widehat F_n(x) - F(x; \widehat \eta_n,1)|, \\
   &  \text{(Anderson-Darling)}  && \mbox{AD} =  \int_{-\infty}^\infty \frac{(\widehat F_n(x) - F(x; \widehat \eta_n,1))^2}{F(x; \widehat \eta_n,1)(1-F(x; \widehat \eta_n,1))} dF(x; \widehat \eta_n,1) 
    \\
&    \text{(Cram\'er-von Mises)}  && \mbox{CM} =  \int_{-\infty}^\infty (\widehat F_n(x) - F(x; \widehat \eta_n,1))^2 dF(x; \widehat \eta_n,1) 
    \\
  &  \text{(Watson)}  && \mbox{WA} = \mbox{CM} -n\Big(\frac1n\sum_{j=1}^n F(Y_j; \widehat \eta_n,1) - \frac12\Big)^2.
\end{align*}
All of these competitor tests were conducted based on the same parametric bootstrap procedure as described in Section~\ref{sec:lnd}.
For this we used the maximum likelihood parameter estimators; see Appendix~\ref{app:pilot_est_scale} for technical details about the practical implementation, also about cases when no maximiser could be found.

\begin{figure}[t]
    \centering
    \includegraphics[width=0.6\textwidth]{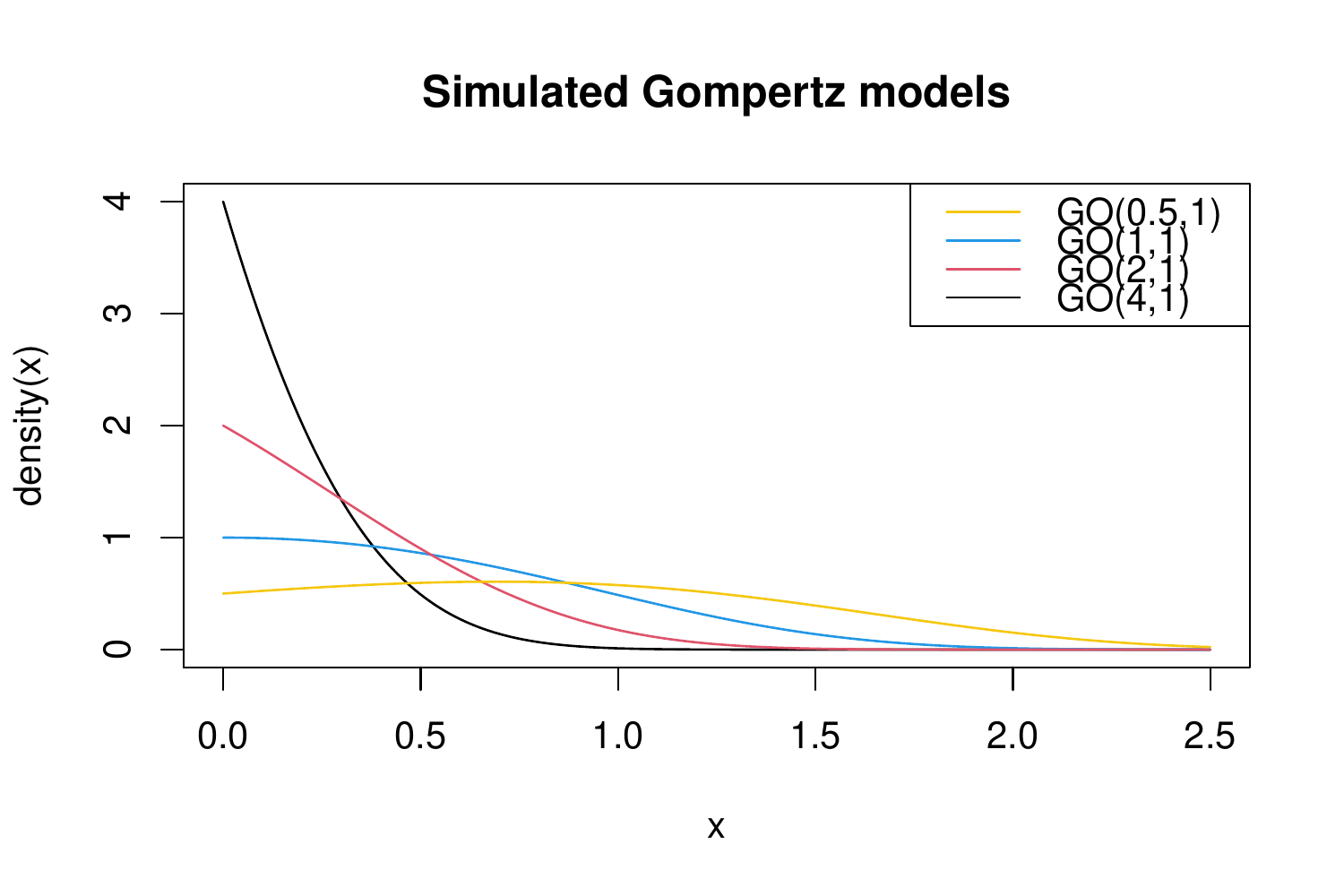}
    
    \caption{Density functions of Gompertz distributions with scale parameter $b=1$ and shape parameters $\eta =$ 0.5 (yellow), 1 (blue), 2 (red), 4 (black).}
    \label{fig:GO-densities}
\end{figure}
The sizes of the tests under the null hypothesis were simulated based on various underlying Gompertz distributions, $GO(0.5,1),GO(1,1), GO(2,1), GO(4,1)$; recall that the tests are scale-invariant, which is why we only let the shape parameter vary. 
Figure~\ref{fig:GO-densities}all illustrates the considered Gompertz distributions in terms of their densities.
In order to compare the tests' power behaviour under the alternative hypothesis, we simulated data according to the non-Gompertz distributions with non-negative support summarised in Table~\ref{tab:alternatives}.
\begin{table}[t]
\centering
\begin{tabular}{|l|l|}
\hline
distribution & density (in $x>0$) \\ \hline
    lognormal $LN(\sigma),\sigma=0.5,1,$ & ${\exp(-\frac{(\log x)^2}{2\sigma^2})}/(\sqrt{2\pi} \sigma x)$, \\[0.2cm]
    Gamma $\Gamma(k) ,k=1,2,3$ & ${x^{k-1} \exp(-x)}/{\Gamma(k)}$, \\[0.2cm]
    inverse Gaussian $IG(\mu,\lambda), \mu =1, \lambda=1,3$ & $\sqrt{\frac{\lambda}{2\pi x^3}} \exp\Big(-\frac{\lambda(x-\mu)^2}{2\mu^2 x}\Big)$, \\[0.2cm]
    Weibull $W(k),k=0.5,3$ & $\frac k\lambda \Big(\frac x \lambda\Big)^{k-1} \exp\Big(-\Big(\frac x\lambda\Big)^k\Big)$, \\[0.2cm]
    Uniform $U(0,5)$ & $0.2 \cdot 1_{(0,5)}(x)$, \\[0.2cm]
    Power $Pow(\nu), \nu=1,2,4$ & $x^{\tfrac1\nu -1} / \nu \ \cdot \  1\{0\leq x \leq 1\}$,\\[0.2cm]
    shifted Pareto $SP(\nu), 3,5,10$  & $\nu (x+1)^{-\nu-1}$,\\[0.2cm]
    linearly increasing failure rate  $LF(\nu), \nu=3,5,10$ & $\nu(x+1) \exp\Big(\frac{1-\nu^2(x+1)^2}{2\nu}\Big)$, \\[0.2cm]
    $GO(1)$-$\Gamma(5)$-mixture $Mix(p), p=0.1,0.25, 0.5, 0.75$ & $p \cdot \exp(1+x- e^{x}) + (1-p)\cdot {x^{4} \exp(-x)}/{\Gamma(5)}$. 
    \\ \hline 
\end{tabular}\\[0.1cm]
\caption{Non-Gompertz distributions considered in the simulation study.}
\label{tab:alternatives}
\end{table}
Note that $\Gamma(1)$ is an exponential distribution ($Exp$),  $Pow(1)$ is a uniform distribution, and, indeed, $LF(a)$ has the linearly increasing failure (or hazard) rate $x\mapsto a(x+1)$.

\begin{table}[!t]
 \centering
 \begin{tabular}{|ll|llllllllll|llll|ll|}
   \hline
  & true & \multicolumn{10}{|c|}{proposed Goodness-of-fit test; tuning parameter $a= $} & \multicolumn{4}{|c|}{classical tests} & \multicolumn{2}{|c|}{not found}  \\
  n & distr. & $0.1$ & $0.25$ & $0.5$ & $0.75$ & $1$ & $1.5$ & $2$ & $3$ & $5$ & $10$ & AD & KS & CM & WA & $\widehat b_n$ & $\widehat b_n^* $ \\ \hline
     20 & $GO(0.5,1)$ & 2 & 3 & 5 & 6 & 6 & 6 & 6 & 6 & 5 & 5 & 5 & 5 & 5 & 5 & 1 & 1 \\ 
   50 &  & 4 & 4 & 5 & 5 & 6 & 6 & 6 & 5 & 5 & 4 & 5 & 5 & 5 & 5 & 0 & 0 \\ 
   100 &  & 4 & 4 & 5 & 5 & 5 & 5 & 5 & 5 & 5 & 5 & 5 & 5 & 5 & 5 & 0 & 0 \\ 
   \hline
   20 & $GO(1,1)$ & 2 & 3 & 4 & 5 & 5 & 6 & 6 & 6 & 6 & 6 & 5 & 5 & 5 & 5 & 3 & 3 \\ 
   50 &  & 3 & 4 & 4 & 4 & 5 & 5 & 6 & 6 & 5 & 5 & 5 & 5 & 5 & 5 & 1 & 2 \\ 
   100 &  & 4 & 4 & 5 & 5 & 5 & 5 & 6 & 6 & 6 & 5 & 5 & 5 & 5 & 5 & 0 & 1 \\ 
   \hline
   20 & $GO(2,1)$ & 2 & 3 & 3 & 4 & 4 & 5 & 6 & 6 & 6 & 6 & 4 & 4 & 4 & 5 & 6 & 7 \\ 
   50 &  & 3 & 3 & 3 & 3 & 4 & 4 & 5 & 5 & 6 & 5 & 5 & 5 & 5 & 5 & 5 & 6 \\ 
   100 &  & 3 & 3 & 3 & 4 & 4 & 5 & 5 & 5 & 6 & 6 & 5 & 5 & 5 & 5 & 2 & 3 \\ 
   \hline
   20 & $GO(4,1)$ & 3 & 3 & 4 & 4 & 4 & 5 & 6 & 6 & 6 & 6 & 4 & 4 & 4 & 4 & 11 & 10 \\ 
   50 &  & 3 & 3 & 3 & 3 & 3 & 4 & 4 & 5 & 6 & 6 & 4 & 4 & 4 & 5 & 10 & 12 \\ 
   100 &  & 3 & 3 & 3 & 3 & 3 & 4 & 4 & 5 & 6 & 6 & 5 & 4 & 4 & 5 & 7 & 10 \\ 
   \hline
 \end{tabular}
 \\[0.2cm]
 \caption{Simulated rounded rejection probabilities (in \%) under the null hypothesis. The considered classical tests are Anderson-Darling (AD), Komogorov-Smirnov (KS), Cram\'er-von Mises (CM), Watson (WA). The last two columns contain the frequency (in \%) of how often the maximum likelihood estimators $\widehat b_n$ and the parametric bootstrap counterpart $\widehat b_n^*$ could not be found.}
 \label{tab:H0}
%
%

\vspace{0.3cm}

 \begin{tabular}{|ll|llllllllll|llll|ll|}
   \hline
  & true & \multicolumn{10}{|c|}{proposed Goodness-of-fit test; tuning parameter $a= $} & \multicolumn{4}{|c|}{classical tests} & \multicolumn{2}{|c|}{not found} \\
  n & distr. & $0.1$ & $0.25$ & $0.5$ & $0.75$ & $1$ & $1.5$ & $2$ & $3$ & $5$ & $10$ & AD & KS & CM & WA & $\widehat b_n$ & $\widehat b_n^*$ \\ \hline
    20 & $LN(0.5)$ & 40 & 47 & 54 & 58 & 60 & 58 & 54 & 45 & 33 & 18 & 53 & 45 & 53 & 54 & 1 & 1 \\ 
   50 &  & 95 & 96 & 97 & 98 & 98 & 98 & 98 & 97 & 93 & 73 & 97 & 93 & 95 & 95 & 0 & 0 \\ 
   100 &  & 100 & 100 & 100 & 100 & 100 & 100 & 100 & 100 & 100 & 100 & 100 & 100 & 100 & 100 & 0 & 0 \\ 
   \hline
   20 & $LN(1)$ & 9 & 10 & 10 & 10 & 11 & 12 & 11 & 11 & 11 & 10 & 15 & 17 & 19 & 19 & 22 & 10 \\ 
   50 &  & 13 & 12 & 14 & 14 & 15 & 15 & 15 & 15 & 15 & 15 & 41 & 32 & 40 & 44 & 34 & 15 \\ 
   100 &  & 15 & 15 & 16 & 18 & 18 & 17 & 17 & 17 & 17 & 17 & 78 & 55 & 67 & 75 & 42 & 21 \\ 
   \hline
   20 & $\Gamma(1)$ & 5 & 6 & 5 & 6 & 6 & 7 & 7 & 7 & 7 & 6 & 7 & 7 & 7 & 6 & 15 & 10 \\ 
   50 &  & 5 & 5 & 6 & 6 & 6 & 6 & 6 & 7 & 6 & 5 & 7 & 7 & 6 & 6 & 19 & 12 \\ 
   100 &  & 6 & 5 & 5 & 6 & 6 & 6 & 6 & 6 & 5 & 5 & 6 & 6 & 6 & 6 & 24 & 16 \\ 
   \hline
   20 & $\Gamma(2)$ & 4 & 5 & 7 & 10 & 11 & 14 & 15 & 16 & 14 & 11 & 11 & 12 & 14 & 15 & 1 & 2 \\ 
   50 &  & 25 & 27 & 31 & 34 & 37 & 41 & 43 & 45 & 46 & 45 & 42 & 32 & 39 & 39 & 0 & 0 \\ 
   100 &  & 61 & 62 & 65 & 68 & 71 & 74 & 76 & 79 & 82 & 84 & 79 & 63 & 71 & 70 & 0 & 1 \\ 
   \hline
   20 & $\Gamma(3)$ & 15 & 18 & 23 & 27 & 29 & 30 & 28 & 24 & 18 & 10 & 23 & 21 & 26 & 26 & 0 & 0 \\ 
   50 &  & 62 & 65 & 69 & 72 & 74 & 75 & 76 & 75 & 72 & 58 & 71 & 57 & 65 & 65 & 0 & 0 \\ 
   100 &  & 94 & 94 & 95 & 96 & 97 & 98 & 98 & 98 & 98 & 97 & 97 & 90 & 94 & 93 & 0 & 0 \\ 
   \hline
   20 & $IG(1,1)$ & 6 & 7 & 8 & 9 & 10 & 12 & 13 & 15 & 17 & 18 & 18 & 18 & 22 & 26 & 14 & 9 \\ 
   50 &  & 9 & 10 & 12 & 15 & 16 & 20 & 24 & 29 & 36 & 45 & 73 & 57 & 62 & 70 & 16 & 11 \\ 
   100 &  & 19 & 17 & 19 & 22 & 24 & 29 & 32 & 38 & 45 & 51 & 99 & 96 & 96 & 98 & 20 & 15 \\ 
   \hline
   20 & $IG(1,3)$ & 31 & 37 & 46 & 52 & 55 & 57 & 55 & 49 & 38 & 23 & 50 & 44 & 51 & 52 & 2 & 1 \\ 
   50 &  & 93 & 94 & 96 & 97 & 97 & 97 & 98 & 97 & 96 & 86 & 98 & 93 & 95 & 95 & 1 & 1 \\ 
   100 &  & 100 & 100 & 100 & 100 & 100 & 100 & 100 & 100 & 100 & 100 & 100 & 100 & 100 & 100 & 0 & 0 \\ 
   \hline
   20 & $W(0.5)$ & 19 & 21 & 21 & 22 & 21 & 23 & 24 & 24 & 25 & 26 & 97 & 91 & 93 & 84 & 51 & 15 \\ 
   50 &  & 20 & 21 & 20 & 21 & 23 & 25 & 28 & 35 & 44 & 54 & 100 & 100 & 100 & 100 & 53 & 19 \\ 
   100 &  & 19 & 20 & 20 & 21 & 23 & 26 & 29 & 39 & 69 & 74 & 100 & 100 & 100 & 100 & 53 & 23 \\ 
   \hline
   20 & $W(3)$ & 15 & 17 & 19 & 17 & 14 & 7 & 3 & 1 & 0 & 0 & 14 & 12 & 15 & 15 & 0 & 0 \\ 
   50 &  & 47 & 48 & 49 & 48 & 45 & 34 & 21 & 9 & 2 & 0 & 40 & 26 & 35 & 35 & 0 & 0 \\ 
   100 &  & 80 & 81 & 82 & 83 & 82 & 77 & 69 & 52 & 26 & 1 & 74 & 50 & 63 & 62 & 0 & 0 \\ 
   \hline
 \end{tabular}
 \\[0.2cm]
 \caption{Simulated rounded rejection probabilities (in \%) under alternative hypotheses with distributions: lognormal ($LN$), gamma ($\Gamma$), inverse Gauss ($IG$), Weibull ($W$), uniform ($U$). The considered classical tests are Anderson-Darling (AD), Komogorov-Smirnov (KS), Cram\'er-von Mises (CM), Watson (WA). The last two columns contain the frequency (in \%) of how often the maximum likelihood estimator $\widehat b_n$ and the parametric bootstrap counterpart $\widehat b_n^*$ could not be found.}
 \label{tab:Ha.1}
 \end{table} 

The results of the simulation study are shown in Tables~\ref{tab:H0}--\ref{tab:Ha.2}. The last two columns therein  indicate how often the maximum likelihood estimator $\widehat b_n$ or its bootstrap counterpart $\widehat b_n^*$ could not be found.
This happened quite rarely under the null hypothesis, with a higher chance for larger shape parameters of the Gompertz distribution (up to 12\% of the iterations for $\eta=4$). 
Under non-Gompertz distributions, these percentages strongly vary from case to case, even within the same family of distributions:
e.g., for Weibull distributions from not at all (Weibull parameter equal to 3) to about 53\% (for $\widehat b_n$, when the parameter equaled 0.5).

Table~\ref{tab:H0} displays the results in terms of empirical rejection rates under the null hypothesis.
We observed only little variation with a change of sample sizes; most of the proposed tests for $a\geq 1.5$ and all of the classical tests exhibited rejection rates very close to 5\%.
However, for tuning parameters $a \leq 1$, most of the proposed tests tend to be conservative with rejection rates going down to 3\%, in some few cases even 2\%.

\begin{table}[t]
 \centering
 \begin{tabular}{|ll|llllllllll|llll|ll|}
   \hline
  & true & \multicolumn{10}{|c|}{proposed Goodness-of-fit test; tuning parameter $a= $} & \multicolumn{4}{|c|}{classical tests} & \multicolumn{2}{|c|}{not found} \\
  n & distr. & $0.1$ & $0.25$ & $0.5$ & $0.75$ & $1$ & $1.5$ & $2$ & $3$ & $5$ & $10$ & AD & KS & CM & WA & $\widehat b_n$ & $\widehat b_n^*$ \\ \hline
  20 & $Pow(1)$ & 9 & 11 & 14 & 15 & 15 & 14 & 13 & 12 & 11 & 10 & 14 & 9 & 11 & 11 & 0 & 1 \\ 
   50 &  & 35 & 36 & 36 & 36 & 33 & 29 & 25 & 21 & 19 & 16 & 35 & 21 & 28 & 28 & 0 & 0 \\ 
   100 &  & 73 & 72 & 70 & 66 & 63 & 55 & 48 & 41 & 34 & 26 & 69 & 43 & 58 & 58 & 0 & 0 \\ 
   \hline
   20 & $Pow(2)$ & 6 & 7 & 9 & 12 & 15 & 20 & 23 & 28 & 34 & 40 & 57 & 30 & 35 & 37 & 12 & 11 \\ 
   50 &  & 21 & 20 & 23 & 29 & 33 & 42 & 48 & 57 & 64 & 71 & 89 & 61 & 71 & 74 & 9 & 11 \\ 
   100 &  & 61 & 53 & 55 & 62 & 66 & 73 & 78 & 83 & 86 & 88 & 99 & 91 & 96 & 97 & 6 & 9 \\ 
   \hline
   20 & $Pow(4)$ & 16 & 15 & 14 & 15 & 15 & 17 & 17 & 17 & 15 & 11 & 99 & 93 & 94 & 92 & 49 & 20 \\ 
   50 &  & 18 & 16 & 14 & 15 & 15 & 16 & 17 & 18 & 17 & 10 & 100 & 100 & 100 & 100 & 52 & 23 \\ 
   100 &  & 20 & 16 & 15 & 16 & 16 & 19 & 20 & 21 & 19 & 13 & 100 & 100 & 100 & 100 & 52 & 25 \\ 
   \hline
   20 & $SP(3)$ & 12 &  12 &  12 & 13 & 13 & 14 & 14 & 14 & 11 & 8 & 31  & 29 & 31 & 20 & 36 & 16 \\
   50 &  & 16 & 16 & 16 & 18 & 19 & 19 & 20 & 19 & 17 & 11 & 54 & 50 & 55 & 35 & 46 & 20 \\ 
   100 &  & 20 & 18 & 19 & 21 & 23 & 24 & 24 & 24 & 21 & 16 & 77 & 72 & 78 & 56 & 51 & 23 \\ 
   \hline
   20 & $SP(5)$ & 9 & 9 & 9 & 9 & 10 & 10 & 10 & 10 & 9 & 7 & 17 & 16 & 18 & 11 & 31 & 17 \\ 
   50 &  & 12 & 11 & 12 & 12 & 13 & 13 & 13 & 13 & 11 & 7 & 28 & 25 & 29 & 16 & 40 & 22 \\ 
   100 &  & 15 & 13 & 14 & 15 & 16 & 17 & 16 & 16 & 14 & 9 & 43 & 39 & 44 & 25 & 47 & 24 \\ 
   \hline
   20 & $SP(10)$ & 7 & 7 & 8 & 8 & 7 & 8 & 8 & 8 & 7 & 6 & 9 & 10 & 10 & 7 & 25 & 16 \\ 
   50 &  & 8 & 8 & 8 & 8 & 8 & 8 & 9 & 9 & 8 & 6 & 13 & 12 & 13 & 9 & 33 & 22 \\ 
   100 &  & 10 & 9 & 9 & 9 & 10 & 10 & 10 & 10 & 9 & 6 & 17 & 15 & 18 & 11 & 38 & 24 \\ 
   \hline
   20 & $LF(3)$ & 3 & 3 & 4 & 5 & 6 & 7 & 7 & 7 & 6 & 5 & 5 & 6 & 6 & 6 & 3 & 4 \\ 
   50 &  & 5 & 5 & 6 & 7 & 8 & 9 & 9 & 9 & 8 & 6 & 7 & 7 & 8 & 8 & 2 & 2 \\ 
   100 &  & 11 & 11 & 12 & 13 & 13 & 14 & 13 & 13 & 12 & 10 & 11 & 10 & 11 & 12 & 0 & 1 \\ 
   \hline
   20 & $LF(5)$ & 3 & 4 & 5 & 5 & 6 & 7 & 7 & 6 & 5 & 4 & 5 & 6 & 6 & 7 & 2 & 3 \\ 
   50 &  & 7 & 8 & 9 & 10 & 10 & 11 & 11 & 11 & 9 & 7 & 8 & 8 & 10 & 10 & 1 & 2 \\ 
   100 &  & 16 & 17 & 17 & 18 & 18 & 19 & 19 & 18 & 16 & 13 & 16 & 14 & 16 & 16 & 0 & 0 \\ 
   \hline
   20 & $LF(10)$ & 4 & 4 & 6 & 7 & 8 & 8 & 8 & 7 & 5 & 3 & 6 & 7 & 7 & 8 & 2 & 2 \\ 
   50 &  & 10 & 11 & 12 & 14 & 14 & 15 & 14 & 13 & 11 & 7 & 12 & 11 & 13 & 13 & 0 & 1 \\ 
   100 &  & 24 & 25 & 26 & 27 & 27 & 28 & 27 & 26 & 24 & 18 & 24 & 19 & 23 & 23 & 0 & 0 \\ 
   \hline
   20 & $Mix(0.1)$ & 10 & 12 & 15 & 16 & 16 & 15 & 13 & 10 & 8 & 7 & 14 & 13 & 15 & 15 & 0 & 0 \\ 
   50 &  & 35 & 36 & 37 & 38 & 38 & 36 & 32 & 26 & 16 & 9 & 34 & 32 & 37 & 36 & 0 & 0 \\ 
   100 &  & 64 & 65 & 65 & 66 & 66 & 64 & 60 & 51 & 33 & 13 & 62 & 59 & 65 & 64 & 0 & 0 \\ 
   \hline
   20 & $Mix(0.25)$ & 5 & 7 & 10 & 13 & 16 & 19 & 21 & 23 & 25 & 24 & 21 & 15 & 15 & 16 & 0 & 1 \\ 
   50 &  & 16 & 19 & 23 & 28 & 32 & 37 & 40 & 44 & 46 & 44 & 40 & 32 & 33 & 33 & 0 & 0 \\ 
   100 &  & 40 & 44 & 50 & 56 & 61 & 67 & 70 & 73 & 73 & 68 & 69 & 59 & 63 & 64 & 0 & 0 \\ 
   \hline
   20 & $Mix(0.5)$ & 9 & 11 & 13 & 15 & 18 & 22 & 24 & 27 & 30 & 32 & 45 & 44 & 48 & 49 & 19 & 9 \\ 
   50 &  & 20 & 20 & 23 & 27 & 30 & 34 & 38 & 43 & 48 & 53 & 85 & 84 & 89 & 90 & 19 & 11 \\ 
   100 &  & 34 & 26 & 28 & 33 & 36 & 42 & 47 & 53 & 59 & 63 & 99 & 99 & 100 & 100 & 17 & 14 \\ 
   \hline
   20 & $Mix(0.75)$ & 18 & 17 & 18 & 19 & 19 & 19 & 18 & 17 & 13 & 9 & 56 & 59 & 64 & 56 & 47 & 15 \\ 
   50 &  & 21 & 19 & 20 & 23 & 23 & 24 & 23 & 22 & 17 & 9 & 95 & 95 & 97 & 95 & 52 & 19 \\ 
   100 & & 21 & 18 & 18 & 23 & 23 & 24 & 23 & 22 & 19 & 16 & 100 & 100 & 100 & 100 & 53 & 23   \\ 
   \hline

 \end{tabular}
 \\[0.1cm]
 \caption{Simulated rounded rejection probabilities (in \%) under alternative hypotheses with distributions: power ($Pow$), shifted Pareto ($SP$), linear failure ($LF$), mixture ($Mix$). The considered classical tests are Anderson-Darling (AD), Komogorov-Smirnov (KS), Cram\'er-von Mises (CM), Watson (WA). The last two columns contain the frequency (in \%) of how often the maximum likelihood estimator $\widehat b_n$ and the parametric bootstrap counterpart $\widehat b_n^*$ could not be found.}
 \label{tab:Ha.2}
 \end{table}

Let us now compare the power results of the proposed and the classical tests; the simulation results are shown in Tables~\ref{tab:Ha.1} and~\ref{tab:Ha.2}.
We generally noticed that the proposed tests are in most cases good competitors of the classical tests in many scenarios; 
at least for some choices of $a$ their rejection rates were close to or slightly greater than those of the classical tests.
This concerns the following distributions: $LN(0.5)$, $\Gamma(k), k=1,2,3$, $IG(1,3)$, and $W(3)$ from Table~\ref{tab:Ha.1} as well as the uniform distribution $Pow(1)$, $LF(\nu), \nu=3,5,10$, and the mixture distributions $Mix(p), p=0.1, 0.25$ from Table~\ref{tab:Ha.2}.
In contrast to that, we also observed that the classical tests clearly outperform the proposed tests in some of the remaining scenarios, e.g.\ for the underlying distributions $LN(1)$, $IG(1,1)$, $W(0.5)$ from Table~\ref{tab:Ha.1} and $Pow(2), Pow(4)$, $SP(3), SP(5)$, $Mix(0.5)$, and $Mix(0.75)$ from Table~\ref{tab:Ha.2}.
It seems that all of the cases go hand in hand with a high chance that the maximum likelihood estimator could not be found.

Most choices of the tuning parameter $a$ resulted in a similar simulated power of the proposed tests; notable exceptions from this can be found for the distributions $LN(0.5)$ (for $n=20$), $\Gamma(2)$, $\Gamma(3)$, $IG(1,1)$, $W(k)$, $Pow(1)$, $Pow(2)$, $Mix(p)$.
The classical tests exhibited a similar power in most scenarios.

All in all, our proposed tests often perform well compared to the classical tests -- but a good choice of the tuning parameter $a$ is of the essence under some alternatives.
Also, if the maximum likelihood estimator for the scale parameter $b$ could not be computed with a relatively high probability, the proposed tests performed suboptimal.

\section{Real data example}\label{sec:real}

In this section, we apply the goodness-of-fit tests to  two different data sets related to lifetimes related to fruitflies and to females born in Germany in 1948 (generated from a life table).
We chose the significance level $\alpha = 5\%$ for all conducted tests.

\subsection{Lifetimes of fruitflies}

\begin{figure}[!t]
    \centering
    \includegraphics[width=1.0\textwidth]{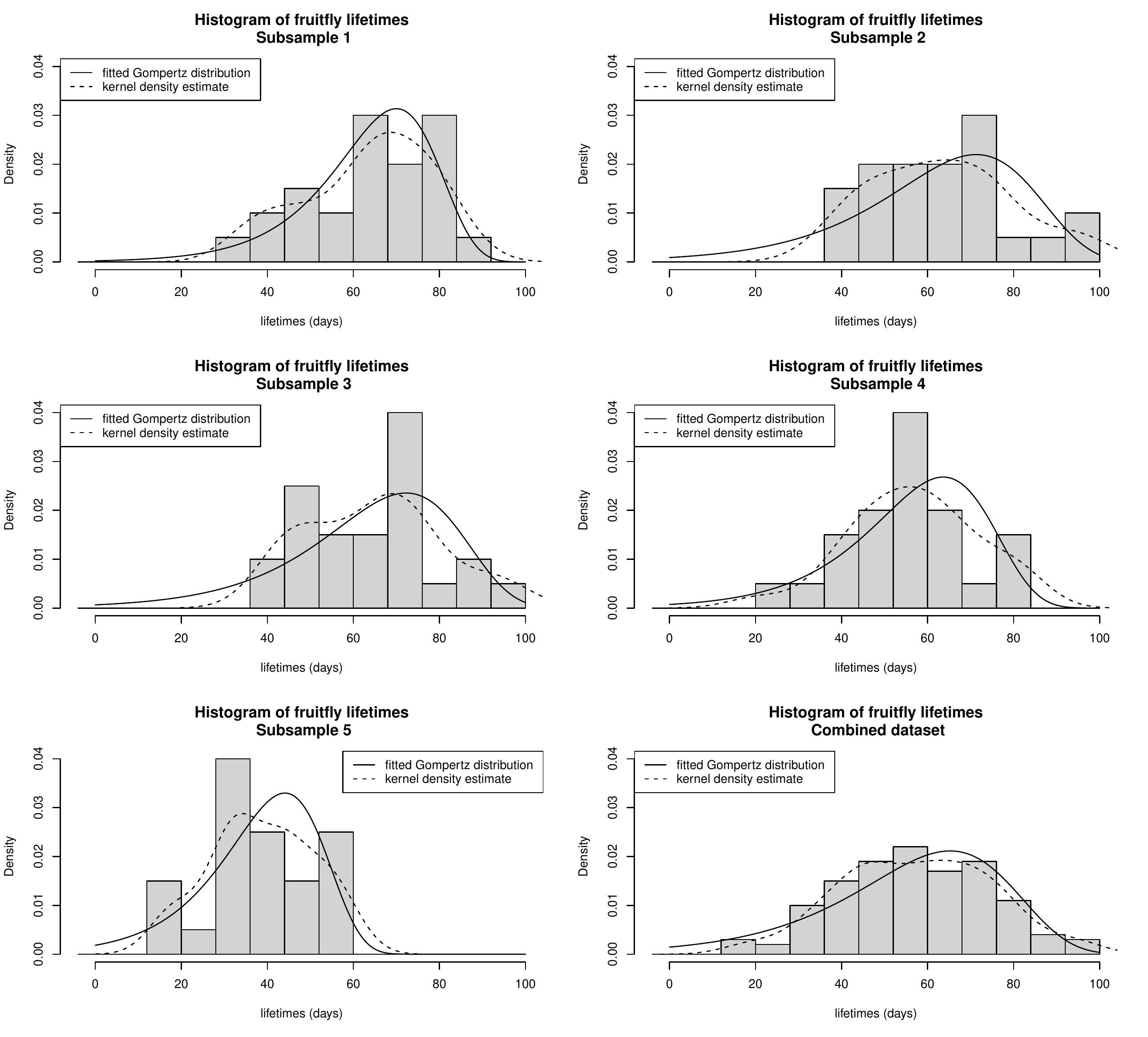}
    \caption{Lifespan (in days) of the fruitflies with separate plots for all subsamples and the complete data set; comparison of histograms, kernel density estimators (default choices in R), and fitted Gompertz distributions.}
    \label{fig:fruitflies_complete}
\end{figure}

The first data set we are about to analyse consists of 125 recorded lifetimes of male fruitflies -- a data set which is publicly available.\footnote{http://jse.amstat.org/jse\_data\_archive.htm; last accessed on August 8, 2022.}
The data resulted from five different groups of male fruitflies of size 25 each; in the different experimental groups varying types and numbers of mating partners.
\cite{PF:1981} argued that ``increasing sexual activity reduces longevity in the male fruitfly''.
For some educational aspects of the data and their analysis, we refer to \cite{H:1983, HS:1994}.

At first, we considered the longevity values of the complete data set, i.e.\ including the measurements from all of the five subgroups.
Figure~\ref{fig:fruitflies_complete} (panel on the bottom-right) illustrates the data set in a histogram which is compared to a nonparametric kernel density estimator and the fitted parametric Gompertz distribution.
The histogram and the kernel density estimator suggest a distribution with at least two modes rather than a unimodal distribution such as the fitted Gompertz distribution.
One reason for the multimodality could be the heterogeneity of the lifetime values in the different groups:
in the other panels of the same Figure, we can see that the location parameters of Subsamples~4 and~5 differ from those of Subsamples~1 to~3.

The $p$-values of the applied goodness-of-fit tests based on 2,000 resampling iterations can be found in Table~\ref{tab:fruitflies}.
It is apparent that, for the complete data set, the proposed test rejects the null hypothesis of Gompertz distributed data for most values of the tuning parameter $a$ -- to be more precise: for all $a\leq 2$.
Also, nearly all of the classical tests arrive at the same result.
Apart from this, almost none of the tests produced a significant outcome when applied separately to the subsamples.
These results could have the following reasons: first, it is not surprising that the power of the tests increase together with the sample size, and the combined data set ($n=125$) is much larger than each of the subsamples ($n=25$).
Second, it is possible that the combination of all five subsamples with most likely different underlying distributions resulted in a sample which cannot be appropriately described by any Gompertz distribution: the mixture of different Gompertz distributions is not a Gompertz distribution.
Also comparing the Gompertz distribution fitted to the combined data set with the kernel density estimator clearly reveals a discrepancy: the density estimate looks almost symmetric and bimodal whereas the fitted Gompertz distribution is left-skewed and unimodal.

As a concluding remark regarding the analysis of these data, we would like to point out that the proposed tests did not produce many surprises when compared to the classical tests.

\begin{table}[t]
 \centering
 \begin{tabular}{|r|rrrrrrrrrr|rrrr|}
   \hline
   &  \multicolumn{10}{c}{proposed Goodness-of-fit test; tuning parameter $a= $} & \multicolumn{4}{|c|}{classical tests} \\ 
   (sub)sample &   $0.1$ & $0.25$ & $0.5$ & $0.75$ & $1$ & $1.5$ & $2$ & $3$ & $5$ & $10$ & AD & KS & CM & WA \\ \hline
   complete & $<\!\textbf{1}$ & $<\!\textbf{1}$ & $<\!\textbf{1}$ & $<\!\textbf{1}$ & $<\!\textbf{1}$ & \textbf{2} & \textbf{4} & 10 & 16 & 36 & $<\!\textbf{1}$ & 5 & \textbf{4} & \textbf{4} \\ 
   1 & 71 & 65 & 55 & 49 & 47 & 48 & 51 & 53 & 48 & 44 & 52 & 49 & 50 & 47 \\
   2 & 6 & \textbf{5} & 5 & 7 & 10 & 27 & 44 & 56 & 55 & 51 & 8 & 27 & 15 & 17 \\
   3 & 6 & 6 & 6 & 8 & 13 & 31 & 47 & 57 & 54 & 49 & 9 & 22 & 15 & 17 \\
   4 & 28 & 27 & 30 & 36 & 44 & 57 & 62 & 60 & 54 & 49 & 22 & 14 & 27 & 28 \\
   5 & 42 & 39 & 35 & 35 & 35 & 43 & 49 & 54 & 56 & 53 & 33 & 9 & 28 & 28 \\
   \hline
\end{tabular}
\\[0.1cm]
\caption{$p$-values (in \%) of all tests applied to the fruitflies data set; rounded to full percentages. Significant results (in view of the significance level $\alpha = 5\%$) are printed in bold-type.}
\label{tab:fruitflies}
\end{table}

\subsection{Data generated from a life table for females born in Germany in 1948}
\label{sec:real_2}

The second type of data sets we are going to analyse is based on life tables provided by the Federal Statistical Office of Germany (Statistisches Bundesamt in Wiesbaden, Germany) published on September 29, 2020, which is publicly available.\footnote{https://www.destatis.de/DE/Themen/Gesellschaft-Umwelt/Bevoelkerung/Sterbefaelle-Lebenserwartung/Publikationen/\_publikationen-innen-kohortensterbetafel.html; last accessed on August 8, 2022.}
Based on the instantaneous hazard rates for females born in 1948 in Germany ($q_x$ in the second column on pp.\ 439-440 of the pdf file that includes the life tables), we reconstructed the underlying probability mass function; see Appendix~\ref{app:hazard_to_pmf} for more details.
Based on these, we could generate data with the help of a multinomial distribution where the probability parameters are equal to the just-mentioned probability mass function.
Because of the early peak due to a relatively high infant mortality and because the hazard rates for the ages above 100 years were aggregated, we decided to crop the distribution to the spans of (i) 10 and 99 years and (ii) 40 and 99 years.
Our idea was to check whether the goodness-of-fit tests are able to detect any deviance from the Gompertz distribution family and whether the fit to some Gompertz distribution is reasonably well if the lifetimes are restricted to all deaths between the ages of 40 and 99, similarly as was motivated in Section~\ref{sec:Intro}.

Based on each of these two truncated distributions, we artificially generated data sets of sizes 20, 50, 100, and 1,000.
Table~\ref{tab:lifetable} contains the $p$-values of the conducted tests.
Let us first focus on the lifetimes truncated to 10 to 99 years.
For $n\leq$\ 100,  none of the classical tests rejected $H_0$ but all of them rejected $H_0$ for $n=$\ 1,000.
None of the proposed tests rejected $H_0$ for $n=20,50$ either.
For $n=$\ 100, the outcomes of the proposed tests do not all agree: all tests with $a \leq 0.75$ rejected $H_0$ but no test with $a\geq 1$ rejected $H_0$.
For $n=$\ 1,000, this threshold for the tuning parameter $a$ is shifted to $a\leq 2$ and $a \geq 3$.

Next, for the lifetimes truncated to 40 to 99 years, we found a similar pattern for the proposed tests, although now fewer of them rejected $H_0$:
those based on $a\geq 0.5$ did not reject $H_0$ for $n=100$.
On the other hand, all of the classical tests rejected $H_0$ for that sample size.
For $n=$\ 1,000, we got similar results, except that now all proposed tests with $a$ up to 1 rejected the null hypothesis; for greater $a$ it was not rejected.
The Kolmogorov-Smirnov test is the only one which rejected $H_0$ for all considered sample sizes.

\begin{figure}[t]
    \centering
    \includegraphics[width=1.0\textwidth]{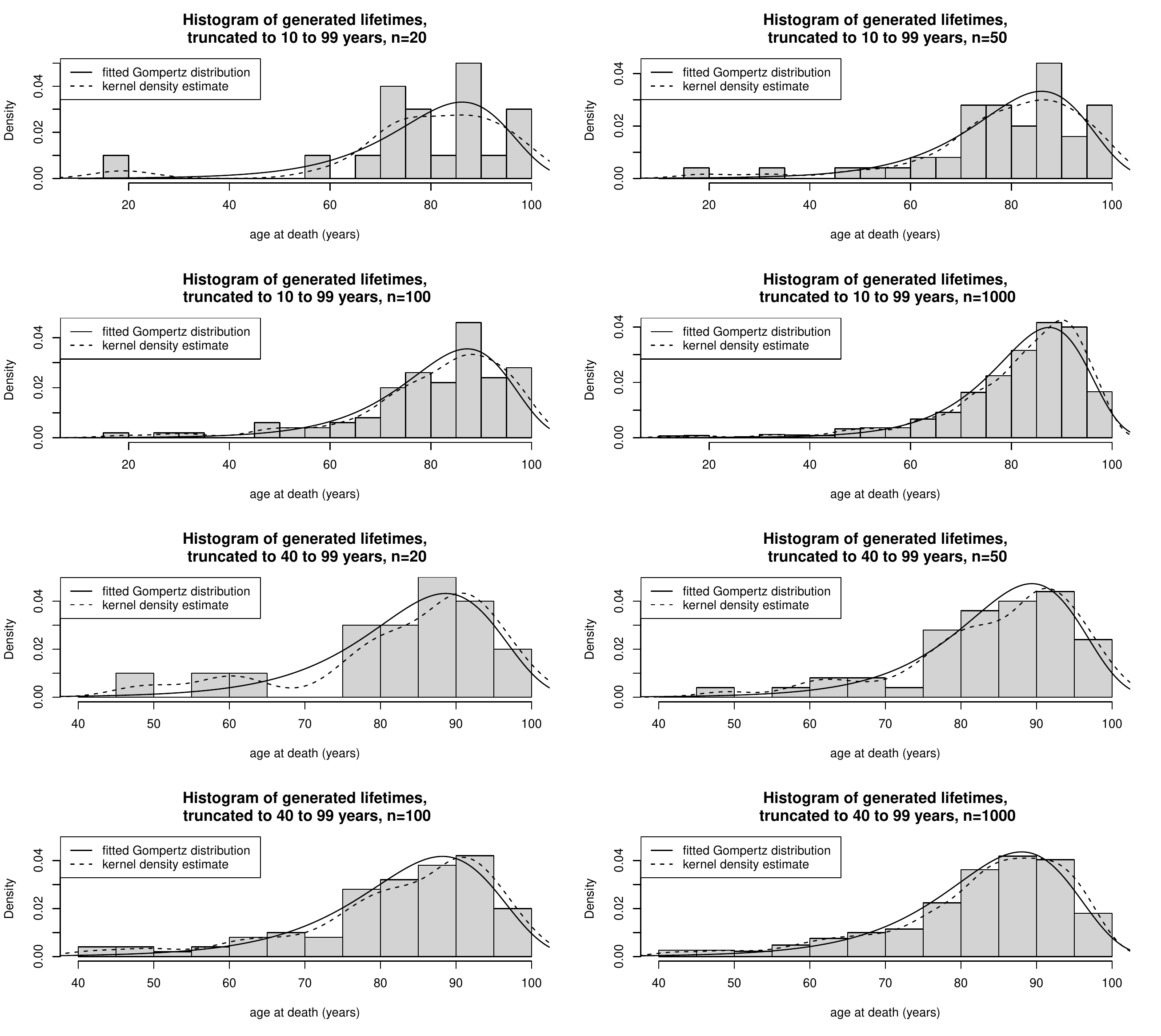}
    \caption{Histograms of artificially generated lifetimes (in years) of women born in Germany in 1948, a kernel density estimator (dashed lines; default values in R), and fitted Gompertz distributions (straight lines); for different truncations of the lifetime distributions: 10 to 99 year (upper half) and 40 to 99 (lower half). The same seed is used to generate all eight data sets.}
    \label{fig:lifetimes_from_tables}
\end{figure}

Some concluding remarks from the perspective of the proposed tests: we see our earlier assumption confirmed by most tests that the lifetimes above 40 are approximately Gompertz-distributed; at least if the tuning parameter $a$ is not overly small.
As the sample size increases to 1,000, more and more tests tend to reject the Gompertz family model which is natural in view of their increasing power and the wrongness of all models.
However, if the truncation is made to the lifetimes between 10 and 99 years, the proposed tests reject $H_0$ more readily, indicating that the Gompertz family model might not be appropriate for the general lifetime distribution;
this remark certainly all the more applies to the completely unrestricted distribution in which a high infant mortality could be observed (not depicted).

\begin{table}[t]
 \centering
 \begin{tabular}{|rr|rrrrrrrrrr|rrrr|}
   \hline
 truncation & sample &  \multicolumn{10}{|c|}{proposed Goodness-of-fit test; tuning parameter $a= $} & \multicolumn{4}{|c|}{classical tests} \\ 
   to years & size &  $0.1$ & $0.25$ & $0.5$ & $0.75$ & $1$ & $1.5$ & $2$ & $3$ & $5$ & $10$ & AD & KS & CM & WA \\ \hline
   10 to 99 & 20 & 18 & 12 & 8 & 7 & 6 & 5 & 5 & 6 & 11 & 39 & 55 & 79 & 74 & 76
   \\
   & 50 &  15 & 12 &  11 &  11 &  11 & 12 &  11 & 12 & 22 & 44 & 55 & 62 & 76 & 82 \\
   & 100 & \textbf{1} & \textbf{1} & \textbf{2} & \textbf{3} & 6 & 10  &  12 & 12 & 26 & 45 & 6 & 7 &  16 & 24 \\ 
    & 1,000 & $<$\textbf{1} & $<$\textbf{1} & $<$\textbf{1} & $<$\textbf{1} & $<$\textbf{1} & $<$\textbf{1} & $<$\textbf{1} & 6 & 9 & 39 & $<$\textbf{1} & $<$\textbf{1} & $<$\textbf{1} & $<$\textbf{1} \\ \hline
    40 to 99 & 20 &  10 & 11 & 14 & 21 & 25 & 28 & 28 & 39 & 39 & 39 & 16 & \textbf{4} & 18 & 19  \\
    & 50 &  16 & 18 & 24 & 30 & 36 & 38 & 37 & 45 & 43 & 43 & 17 & \textbf{3} & 12 & 13 \\
    & 100 &  \textbf{2} & \textbf{3} & 5 & 9 & 15 & 25 &   31 & 47 & 45 & 44 & \textbf{2} & $<$\textbf{1} & \textbf{3} & \textbf{3} \\ 
                & 1,000 & $<$\textbf{1} & $<$\textbf{1} & $<$\textbf{1} & $<$\textbf{1} & $<$\textbf{1} & 12 & 33 & 66 & 58 & 52 & $<$\textbf{1} & $<$\textbf{1} & $<$\textbf{1} & $<$\textbf{1} \\ \hline
\end{tabular}
\\[0.1cm]
\caption{$p$-values (in \%) of all tests applied to the data sets generated based on the (truncated) life tables; rounded to full percentages. Significant results (in view of the significance level $\alpha = 5\%$) are printed in bold-type.}
\label{tab:lifetable}
\end{table}

\section{Discussion and outlook}\label{sec:CaO}

In the present paper, we demonstrated how a Stein characterisation for the Gompertz distribution family can be used to develop a parametric bootstrap-based goodness-of-fit test for the composite null hypothesis.
In our simulation study and the real data analyses, we have used the maximum likelihood estimators of the parameters of the Gompertz distribution.
Other choices of parameter estimators could also be covered by the developed theory as long as they exhibit an asymptotically linear structure, for more information on the influence of parameter estimation techniques on the power of goodness-of-fit tests see \cite{DKO:1990}.
This could potentially solve the issue of the suboptimal power seen in Tables~\ref{tab:Ha.1}
and~\ref{tab:Ha.2} whenever there was a high probability that the maximum likelihood estimator for the scale parameter $b$ could not be found.
Other than that, the proposed test revealed a good control of the type-I error probability across multiple underlying Gompertz distributions and a satisfactory power behaviour under many considered alternative hypotheses, also when compared to classical competitor tests.
Still, the choice of the tuning parameter $a$ is crucial for obtaining a reliable test.
In general, if no further information is available, intermediate choices of $a$ close to $1.5$ seem to be safest.
Another possibility to solve this issue is to combine the proposed test with an adaptive selection procedure; see for instance \cite{T:2019} for a bootstrap-based approach.
It should be noted that the additional bootstrap layer would significantly increase the computational complexity of the test procedure.

Finally, in view of applications to medical time-to-event data, another important extension of the proposed test would involve the handling of censored data.
The difficulties in this connection are two-fold:
firstly, the test statistic would need to involve the Kaplan-Meier estimator of the re-scaled observations instead of the empirical cumulative distribution function and, in particular, the expectation given in the Stein characterisation needed to be replaced by another estimator for which there is no standard approach.
Secondly, the maximum likelihood estimators of the Gompertz distribution parameters would change.
The large sample properties of the resulting test statistic could potentially be established by means of adaptations of techniques from survival analysis.

\section*{Acknowledgements}
The authors thank B.\ Clau{\ss} for preliminary work on this topic in his master thesis and O.\ Ku\ss\  for helpful discussions.

\section*{Conflict of interest statement}
Both authors declare that there are no  financial or commercial conflicts of interest.

\bibliographystyle{abbrv}
\bibliography{lit-GO}

\newpage
\begin{appendix}
\section*{Appendix}

\section{Proofs}\label{app:PT1}
First we show a Stein-type characterisation of the Gompertz distribution, which is a special case of the so called density approach in Stein's method, see \cite{SDHR:2004}, Proposition 1.4. The proof of Lemma \ref{lem:char} follows the lines of proof of Theorem 1 in \cite{BE:2019} and is hence omitted.

\begin{lemma} \label{lem:char}
   A positive random variable $X$ is $GO(\eta, b)$, $\eta,b>0$, distributed, if and only if
   \begin{align*}
       \mathbb{E} \big[g'(X) + \big(-\eta b e^{bX} + b\big)g(X)\big] = 0
   \end{align*}
   holds for all functions $g \in \mathcal{G}$, where
   \begin{align*}
   \mathcal{G} = \big\{g: (0,\infty) \rightarrow \mathbb{R}| &\;g \text{ is differentiable}, g' (x) \text{ and } \big(\eta b e^{bx} - b\big) g(x) \text{ are bounded, and}\\
   &\lim_{x \downarrow 0} g(x) = \lim_{x \downarrow 0} g(x) f(x) = \lim_{x \rightarrow \infty} g(x) f(x) = 0\big\},
   \end{align*}
   and $f$ represents the probability density function of $GO(\eta, b)$ in \eqref{eq:density}.
\end{lemma}
\emph{Proof of Theorem \ref{thm:char}.}  
First assume that $X \sim GO(\eta, b)$. With $f(x)=f(x;\eta,b)$ as in \eqref{eq:density} we have
\begin{align*}
    \mathbb{E}\big[(\eta b e^{bX} - b)\mathbf{1}\{X > s\}\big] &= \int_0^\infty (\eta b e^{bx} - b) \mathbf{1}\{x > s\} f(x)\, dx= - \int_0^\infty f'(x) \mathbf{1}\{x > s\}\,dx = f(s),\quad s>0.
\end{align*}
With $\min\{X, s\} = \int_0^s \mathbf{1}\{X > t\} \,dt$ and the theorem of Fubini, we have
\begin{align*}
    T^X (s) &= \mathbb{E}\big[(\eta b e^{bX} - b)\min\{X, s\}\big] = \int_0^\infty (\eta b e^{bx} - b) \int_0^s \mathbf{1}\{x > t\}\,dt f(x)\,dx\\
    &=\int_0^s \int_0^\infty (\eta b e^{bx} - b) \mathbf{1} \{x > t\} f(x) \,dx \,dt = \int_0^s f(t) \,dt = F(s).
\end{align*}
for all $s > 0$ and hence $T^X \equiv F$ on $\mathbb{R}$.\\
Assume now \(T^X \equiv F\). Define 
$d^X (s) = \mathbb{E}\big[(\eta b e^{bx} - b) \mathbf{1}\{X > s\}\big] \mathbf{1}\{0 < s < \infty\}$ for all \(s \in \mathbb{R}\).
By the theorem of Tonelli, we have
\begin{align*}
    \int_0^\infty \mathbb{E}\big[|\eta b e^{bX} - b| \mathbf{1}\{X > t\}\big]\,dt &= \mathbb{E} \Big[|\eta b e^{bX} - b| \int_0^\infty \mathbf{1}\{X > t\}\,dt\Big]\\
    &= \mathbb{E}\big[|\eta b e^{bx} - b|X\big] \leq \eta b \mathbb{E}\big[Xe^{bX}\big] + b\mathbb{E}\big[X\big] < \infty.
\end{align*}
Hence the theorem of Fubini is applicable and
\begin{align*}
    \int_{- \infty}^s d^X (t) \,dt &= \int_0^s \mathbb{E} \big[(\eta b e^{bX} - b) \mathbf{1}\{X > t\}\big] \,dt= \mathbb{E} \Big[(\eta b e^{bX} - b) \int_0^s \mathbf{1}\{X > t\} \,dt\Big]\\
    &= \mathbb{E} \big[(\eta b e^{bX} - b) \min\{X, s\}\big] = T^X (s),\quad s\in\R,
\end{align*}
follows, and \(d^X \geq 0\) almost everywhere, since \(T^X\) is monotonically increasing, and
\[\int_\mathbb{R} d^X (t)\,dt = \lim_{s \rightarrow \infty} T^X (s) = \lim_{s \rightarrow \infty} F(s) = 1.\]
It follows that \(d^X\) is the probability density function of \(T^X\), and hence \(F\). Let \(g \in \mathcal{G}\) with \(\mathcal{G}\) from Lemma \ref{lem:char}. Since $g$ is bounded, we have by Fubini
\begin{align*}
    \mathbb{E}[g'(X)] &= \int_0^\infty g'(s) d^X(s) \,ds = \int_0^\infty f'(s) \mathbb{E}\big[(\eta b e^{bx} - b) \mathbf{1}\{X > s\}\big] \,ds\\
    &= \mathbb{E}\big[ (\eta b e^{bx} - b) \int_0^\infty g'(s) \mathbf{1}\{X > s\} \,ds\big] = \mathbb{E}\big[ -(-\eta b e^{bx} + b) (g(X) - g(0))\big]\\
    &= \mathbb{E}\big[ -(-\eta b e^{bx} + b) g(X)\Big].
\end{align*}
By Lemma \ref{lem:char} the claim follows.\hfill$\square$\\[2mm]
\emph{Proof of Theorem \ref{thm:H0vert}.} 
The key idea is to use the central limit theorem in Hilbert spaces. Since in \eqref{eq_V} we have a sum of dependent random variables due to the estimators and rescaled data the CLT is not directly applicable. Hence, we first introduce the helping processes
\begin{align*}
    \widetilde{V}_n (s) = &\frac{1}{\sqrt{n}} \sum_{j=1}^n \Big\{ (\eta_n e^{X_{n, j}} - 1) \min \{ X_{n, j}, s\} - \mathbf{1} \{X_{n, j} \leq s\}\Big\} \\
    &+ \sqrt{n} (\widehat{\eta}_n - \eta_n) \frac1n \sum_{j=1}^n \Big\{ e^{X_{n, j}} \min \{X_{n, j}, s\} \Big\}\\
    &+ \sqrt{n} (\widehat{b}_n - 1) \frac{1}{n} \sum_{j=1}^n \Big\{ \big(\eta_n e^{X_{n, j}} + \eta_n X_{n, j} e^{X_{n, j}} - 1\big) \min \{X_{n, j}, s\}\Big\}, \quad s>0,
\end{align*}
and
\begin{align*}
    V_n^* (s) = & \frac{1}{\sqrt{n}} \sum_{j=1}^n \Big\{ (\eta_n e^{X_{n, j}} - 1) \min \{ X_{n, j}, s\} - \mathbf{1} \{X_{n, j} \leq s\} \\
    &+ \psi_1 (X_{n, j}, \eta_n) \mathbb{E} \big[e^X \min \{X, s\}\big] \\
    &+ \psi_2 (X_{n, j}, \eta_n) \Big\{ \eta_n \mathbb{E}\big[e^X \min \{X, s\}\big] + \eta_n \mathbb{E}\big[Xe^X\min \{X, s\}\big] - \mathbb{E}\big[\min \{X, s\}\big] \Big\} \Big\}, \quad s>0.
\end{align*}
By a multivariate Taylor approximation around $(\eta_n,1)$, some integral transform, the use of H\"{o}lders inequality and conditions \eqref{eq:w1} and \eqref{eq:w2}, repeated use of Slutzki's Lemma and tightness arguments, the consistency of the estimators, as well as lengthy calculations, we are able to show that $\|\sqrt{n}V_n-\widetilde{V}_n\|_{\mathbb{H}}=o_{\PP}(1)$ and $\|\widetilde{V}_n-V_n^*\|_{\mathbb{H}}=o_{\PP}(1)$. Denote the $j$-th summand of $V_n^*$ by $W_{n,j}$. Then direct calculation shows $\E[W_{n,j}]=0$ and we have a sequence of rowwise identically distributed random variables. We have
\begin{align} \label{eq:chen3.1b}
    \lim_{n\rightarrow \infty} \mathbb{E} \Big\Vert \frac{1}{\sqrt{n}} \sum_{j=1}^n W_{n, j} \Big\Vert_{\mathbb{H}}^2 < \infty.
\end{align}
The central limit theorem of Lindeberg-Feller shows for all \(g \in \mathbb{H}\setminus \{0\} \) 
\begin{align} \label{eq:chen3.1a}
    \frac{1}{\sqrt{n}} \sum_{j=1}^n \langle W_{n, j}, g \rangle_{\mathbb{H}} \stackrel{\mathcal{D}}{\longrightarrow} \mathcal{N}(0, \sigma_\eta^2 (g)),\quad n\rightarrow\infty,
\end{align}
where $\sigma_\eta^2 (g) = \lim_{n \rightarrow \infty} \mathbb{E}\big[\langle W_{n, 1}, g\rangle_\mathbb{H}^2\big]$.
By Lemma 3.1 and Remark 3.3 in \cite{CW:98}, \eqref{eq:chen3.1b} and \eqref{eq:chen3.1a} we have
\[\frac{1}{\sqrt{n}} \sum_{j=1}^n W_{n, j} \stackrel{\mathcal{D}}{\longrightarrow} \mathcal{W},\]
where \(\mathcal{W}\) denotes a centered Gaussian random element of \(\mathbb{H}\) with covariance operator \(\mathcal{C}\), which satisfies for all \(g \in \mathbb{H} \setminus \{0\}\) the equation \(\sigma_\eta^2(g) = \langle \mathcal{C} g, g \rangle_{\mathbb{H}}\). The covariance operator is identical to $\lim_{n\rightarrow\infty}\mathbb{E}[W_{n, 1} (s) W_{n, 1} (t)]$, which after a considerable amount of straightforward calculation provides the stated formula in the theorem. Note that we derived and used the identities
\begin{eqnarray*}
    \mathbb{E}\big[(\eta e^{X} - 1)^2 \min\{X, s\} \min\{X, t\}\big] &=& -sf(s; \eta,1) - sf(t; \eta,1) + sF(s; \eta,1)+ s(1-t)F(t; \eta,1)\\
    && + \mathbb{E}[X\mathbf{1}\{X \leq s\}] + \mathbb{E}[X^2\mathbf{1}\{X \leq s\}] + \mathbb{E}[X\mathbf{1}\{X \leq t\}],\\
    \mathbb{E}\big[(\eta e^{X} - 1) \min\{X, s\} \mathbf{1}\{X \leq t\}\big] &=& - sf(t; \eta,1) + F(s; \eta,1)\\
\mathbb{E}\big[(\eta e^{X} - 1)X \mathbf{1}\{X \leq s\}\big] &=& -sf(s; \eta,1) + F(s; \eta,1)\\
\mathbb{E} \big[Xe^X \min\{X, s\}\big] &=& \frac{1}{\eta}\big(s - sF(s; \eta,1) + \mathbb{E}\big[X\big] + \mathbb{E}\big[X\mathbf{1}\{X \leq s\}\big]+ \mathbb{E}\big[X^2\mathbf{1}\{X \leq s\}\big]\big)\\
    \mathbb{E} \big[e^X \min\{X, s\}\big] &=& \frac{1}{\eta}\big(1 + \mathbb{E} \big[X \mathbf{1}\{X \leq s\}\big]\big)\\
    \mathbb{E} \big[\min\{X, s\}\big] &=& s - sF(s; \eta,1) + \mathbb{E} \big[X \mathbf{1}\{X \leq s\}\big],
\end{eqnarray*}
where $X\sim GO(\eta,1)$ and $0<s\le t<\infty$. Then, Slutzki's lemma, the continuous mapping theorem combined with the triangular inequality prove the statement. \hfill $\square$.\\[2mm]
%
%
\emph{Proof of Theorem~\ref{satz:boot1}.}
Write $H^{\eta}(\cdot)$ for the distribution of $\|\mathcal{W}\|_{\mathbb{H}}^2$ and $H_n^{\eta_n}$ for the distribution of $T_n$. Note that $H^{\eta}$ is continuous and strictly increasing on $\{s>0: 0< H^{\eta}(s) < 1\}$. By Theorem \ref{thm:H0vert} it holds that $H_n^{\eta_n}(s) \rightarrow H^{\eta}(s)$ for each $s>0$ as $n\rightarrow\infty$, so by continuity of $H^{\eta}$ we have
\begin{equation*}
\sup_{s>0}\left|H_n^{\eta_n}(s)-H^{\eta}(s)\right|\longrightarrow0\quad \mbox{as}\,n\rightarrow\infty.
\end{equation*}
A combination of the last result with the consistency of $\widehat{\eta}_n$ yields
\begin{equation*}
\sup_{s>0}\left|H_n^{\widehat{\eta}_n}(s)-H^{\eta}(s)\right|\cp0\quad \mbox{as}\,n\rightarrow\infty.
\end{equation*}
Hence, we have by an identical construction as in (3.10) of \cite{H:1996}
\begin{equation*}
\sup_{s>0}\left|H_{n, B}^*(s)-H^{\eta}(s)\right|\cp0\quad \mbox{as}\,n,B\rightarrow\infty,
\end{equation*}
from which $c_{n,B}^*(\alpha)\cp \inf\{s:H^{\eta}(s)\ge 1-\alpha\}$ follows as $n,B\rightarrow\infty$ . This implies the claim. \hfill$\square$
\\[0.2cm]
\emph{Proof of Theorem \ref{thm:cons}.}
Denote 
\begin{align*}
    \bar{V}_n (s) = &\frac{1}{n} \sum_{j=1}^n \Big\{ (\eta_0 e^{X_{j}} - 1) \min \{ X_{j}, s\} - \mathbf{1} \{X_{j} \leq s\}\Big\} +  (\widehat{\eta}_n - \eta_0) \frac1n \sum_{j=1}^n \Big\{ e^{X_{j}} \min \{X_{j}, s\} \Big\}\\
    &+  (\widehat{b}_n - 1) \frac{1}{n} \sum_{j=1}^n \Big\{ \big(\eta_0 e^{X_{j}} + \eta_0 X_{j} e^{X_{j}} - 1\big) \min \{X_{j}, s\}\Big\}, \quad s>0,
\end{align*}
By the same reasoning as in the proof of Theorem \ref{thm:H0vert}, we have $\|V_n-\bar{V}_n\|_{\mathbb{H}}=o_{\PP}(1)$. Hence by the triangle inequality, we have
\begin{equation*}
\frac{T_n}{n}=\|V_n\|_{\mathbb{H}}^2\le\|V_n-\bar{V}_n\|_{\mathbb{H}}^2+\|\bar{V}_n\|_{\mathbb{H}}^2= \|\bar{V}_n\|_{\mathbb{H}}^2+o_{\PP}(1)
\leq  \|{V}_n\|_{\mathbb{H}}^2+o_{\PP}(1).
\end{equation*}
Due to the consistency of $\hat{\eta}_n$ and $\hat{b}_n$ for $\eta_0$ and 1, respectively, and the finiteness of all required expectations, the last two terms in the definition of $\bar V_n$ converge to 0 in probability by the law of large numbers and Slutzki's lemma.
This convergence holds uniformly in $s>0$.
Furthermore, due to monotonicity arguments in combination with the law of large numbers, the first term in the definition of $\bar V_n$ converges to $\Delta_{\eta_0}^*(s)$ in probability, uniformly in $s>0$.
It follows that 
\begin{align*}
\|\bar{V}_n\|_{\mathbb{H}}^2
= \int_0^\infty \{ \bar{V}_n^2(s) - \Delta^{*2}_{\eta_0}(s)\} w(s) \textnormal{d}s + \Delta_{\eta_0}.
\end{align*}
Here, the first term is bounded in absoulte value by $\sup_{s>0}|\bar{V}_n^2(s) - \Delta^{*2}_{\eta_0}(s)| \cdot \|1\|_{\mathbb{H}}^2 = o_{\PP}(1)$ due to the uniform convergence in probability argued above, and $\|1\|_{\mathbb{H}}^2<\infty$ due to \eqref{eq:w1}.
Another application of Slutzki's lemma concludes the proof.\hfill $\square$.

\emph{Proof of Corollary~\ref{cor:cons}.}
Due to the assumption given in \eqref{eq:psi_11}--\eqref{eq:psi3}, we know that the estimators $\widehat \eta_n$ and $\widehat b_n$ converge in probability to some $\tilde\eta >0$ and $\tilde b >0$, respectively. 
Hence, for a given level of significance $\alpha$, we know from the proof of Theorem~\ref{satz:boot1} that the critical values $c_{n, B}^* (\alpha)$ converge to a fixed value $0 < c < \infty$ (say) for $n, B \rightarrow \infty$. Since $\Delta_{\eta_0}$ from  Theorem \ref{thm:cons} is strictly positive by the characterisation in Theorem \ref{thm:char} if the underlying law is not from the Gompertz family $\mathbf{\mbox{GO}}$, the claim follows directly.
\hfill$\square$

\section{More details on the maximum likelihood estimation from a practical point of view}
\label{app:pilot_est_scale}

Since there is no closed-form solution for the maximum likelihood estimators $\widehat \eta_n$ and $\widehat b_n$ in the Gompertz family, we have used the Newton-Raphson algorithm to approximate the maximiser numerically; we used the R-package \emph{pracma} for this.
This algorithm requires an initial guess for the scale parameter $b$ in order to find $\widehat b_n$.
Next, the maximum likelihood estimator for $\eta$ depends on $\widehat b_n$ through $\widehat \eta_n = (\frac1n \sum_{j=1}^n \exp(\widehat b_n x_j) -1 )^{-1}$.

In this appendix, we will explain how to find such a pilot estimator.
Our idea was to involve an estimator of the cumulative hazard function, $\Lambda : x \mapsto \eta (\exp(bx) -1 )$.
In particular, it is well-known that 
$$ \widehat \Lambda_n: x \mapsto \sum_{j: X_{(j)} \leq x} \frac{1}{n-j+1} $$ 
is consistent for $\Lambda(x)$; here, $X_{(1)}, \dots, X_{(n)}$ denote the order statistics.
Next, we chose a rather large value $\hat z_n$ of the data $X_1, \dots, X_n$, e.g.\ their 90$^{th}$ percentile; 
in any case, $\hat z_n$ should converge to some fixed and finite value $z>0$ as the sample size increases.
Our pilot estimator is then given as 
$$ \widehat{b}_{n,\text{pilot}} = \frac{2}{\hat z_n} \log \Big( \frac{\widehat \Lambda_n(\hat z_n) - \widehat \Lambda_n(\hat z_n/2)}{\widehat \Lambda_n(\hat z_n/2)}  \Big). $$
Indeed, due to the uniform consistency of $\widehat \Lambda_n$ for $\Lambda$ on compact intervals, $\widehat{b}_{n,\text{pilot}}$ is a consistent estimator for
$$ \frac{2}{z} \log \Big( \frac{ \Lambda(z) -  \Lambda(z/2)}{ \Lambda(z/2)}  \Big)
= \frac2z \log \Big( \frac{\exp(bz) - \exp(bz/2)}{\exp(bz/2) -1} \Big)
= \frac2z \log ( \exp(bz/2) \cdot 1) = b. $$

Thus, the scale parameter estimator was found as the solution to 
\begin{align}
\label{eq:MLE}
h(b; x_1, \dots, x_n) = \Big(\frac1n \sum_{j=1}^n \exp(b x_j) - 1\Big)\cdot (b \bar  x_n + 1) - \frac bn \sum_{j=1}^n  x_j \exp(b x_j) = 0     
\end{align}
with the help of the Newton-Raphson algorithm.
We involved an estimator of the Gompertz cumulative hazard function $x\mapsto \eta (\exp(bx)-1)$ in order to find a reasonable pilot estimator as an initial value for $b$ in the algorithm; details on the pilot estimator can be found in Appendix~\ref{app:pilot_est_scale}.
Whenever \eqref{eq:MLE} had no solution for $b$ in the positive numbers,
we have chosen the rather small value $\widehat b_n = 0.001$, 
since $\lim_{b \downarrow 0} h(b; x_1, \dots, x_n) = 0$.

\section{Reconstruction of a probability mass function from a discrete hazard rate}
\label{app:hazard_to_pmf}

In Section~\ref{sec:real_2}, we generated data according to an official life table. In this appendix, we will explain the procedure in detail.

Let $X$ be a random variable with values in $\mathbb{N}_0$ and probability mass function, $k \mapsto p(k) = P(X=k)$, $k \in \mathbb{N}_0$.
It is well-known that the hazard function is obtained as 
$$ q(k) = P(X=k \ | \ X \geq k) =  \frac{p(k)}{S(k-1)}, $$
where $S(k) = P(X > k) = \sum_{\ell = k+1}^\infty p(\ell) = \prod_{\ell = - \infty}^k (1 - q(\ell))$ is the so-called survival function.
Thus, for a given hazard function $q$, the corresponding probability mass function can be reconstructed based on the following iterative procedure:
\begin{itemize}
    \item $p(0) = q(0)$,
    \item $S(k-1) = \prod_{\ell=0}^{k-1} (1 - q(\ell)) $, \ $k \geq 1$,
    \item $p(k) = S(k-1) \cdot q(k)$, \ $k \geq 1$.
\end{itemize}
For each $k$, the last two steps have to be conducted one after the other, before increasing $k$ to the next integer value.

For the data set in Section~\ref{sec:real_2}, a final adjustment was necessary to ensure that the probability mass function adds up to~1; due to rounding errors this was not immediately the case.
Thus, the final step is to take $\tilde p(k) = p(k)/\sum_{\ell=0}^\infty p(\ell)$, $k\in \mathbb{N}_0$, as the probability mass function.

In Section~\ref{sec:real_2}, we also considered some truncated distributions (with probability mass functions, say, $\check p$) whose probability mass functions were obtained as follows:
for each $k$ within the truncation region, say, $k \in [0, L] \cup [R,\infty)$ with integers $0 \leq L <  R < \infty$, we set $\check p (k) = 0$.
For $k \in [L+1, R-1]$, we set $\check p(k) = \tilde p(k)/\sum_{\ell=L+1}^{R-1} p(\ell)$, to rescale $\check p$ to a probability mass function.

\end{appendix}
\end{document}